\documentclass[11pt,a4paper,reqno, oneside]{article}
\usepackage{latexsym,amsfonts,amsthm,amsmath,amssymb}
\usepackage{amstext}
\usepackage{xspace}
\numberwithin{equation}{section}
\theoremstyle{plain}
\newtheorem{teo}{Theorem}[section]
\newtheorem{lem}[teo]{Lemma}
\newtheorem{prop}[teo]{Proposition}
\newtheorem{cor}[teo]{Corollary}
\newtheorem*{ack}{Acknowledgments}
\theoremstyle{definition}
\newtheorem{defi}[teo]{Definition}
\newtheorem{rem}[teo]{Remark}

\newcommand{\R}{\mathbb{R}}
\newcommand{\N}{\mathbb{N}}

\newcommand{\Om}{\Omega}
\begin{document}
\title{Variational Problems on Classes of Convex Domains}
\author{Nicolas Van Goethem}
\date{}
\maketitle

\begin{abstract}
 We prove the existence of  minimizers for functionals defined over the class  of convex domains contained inside a bounded set $D$ of $\R^N$ and with prescribed volume. Some applications are given, in particular we prove that the  eigenvalues of differential operators of second and fourth order with non-constant coefficients as well as integral functionals depending on the solution of an elliptic equation can be minimized over this class of domains. Another application of this result is related to the famous Newton problem of minimal resistance. In general, all the results we shall develop hold for elliptic operators of any even order larger that $0$.
\end{abstract}
\section{Introduction}
Minimizing functionals over a class of domains is a problem with many applications. A classical application is e.g. Newton's problem of minimal resistance (see \cite{BuKa}, \cite{Ka}) where one seeks an optimal shape of a body $\Om$ in a class of sets, for instance convex domains of $\R^N$. In this case, the natural functional to be associated with $\Om$ is the so called total resistance of the body, i.e. an integral on $\Om$ involving a suitable function of the normal vector at each point of the boundary. The integrand may be related to the solution of a simple PDE (see \cite{Ka} for Newton's problem) but sometimes optimization problems  are considered more general, for instance the minimization on a space of measures (see \cite{BuDa} for the shape optimization problem). Then, particular tools from geometric measure theory can be required to solve the problem. 

This paper is based on some simple concepts related to convexity. One of the concepts we need in an optimization problem is the space on which we minimize the given functional: the larger the class of candidates for the minimizer is, the harder the proof of an existence result is. One is often forced to introduce constraints on the space, and these should be natural enough, otherwise one loses some of the relevance of the problem. For example,  the concept of compact sets is often the key of the direct method of the calculus of variations but since compactness is a very strong condition, how can it be obtained for general spaces? As is often the case, one has to pay for the generality of a result, and the cost is precisely the complexity of the topology on the chosen space. This is a reason why some strong tools from geometric measure theory are used if one seeks to minimize shape functionals on the space of all open sets of $\R^N$ (see \cite{BuDa}). As soon as the domains are general enough, the candidates for the optimality may have very different kinds of behaviors. This variety of candidates can be compensated either by constraining the class of functionals, i.e. only a narrow class of functionals is proved to assume their minima, or by the choice of the ``most suitable'' topology adapted to the problem. Eventually, the chosen topology should give some kind of compactness. Alternatively, one can also focus on the minimizing sequences to prove existence results. 

In a recent work of Kawohl (see \cite{Ka}), this second approach is used for proving the existence of a minimum on a class of convex domains, provided the functional is bounded from below, homogeneous of degree $-s$, monotone decreasing with respect to the set inclusion and coercive. The two latter conditions allow the author to consider domains which are not a priori bounded. In this paper we consider bounded domains since one is often lead anyhow to consider this situation. When considering elliptic equations with non-constant coefficients, the homogeneity and the shift invariance of the functionals are too hard constraints. Moreover, it seems that assuming both convexity of the domains and  monotonicity of the functional is redundant. Indeed, it can be showed that monotonicity can sometimes  be compensated by convexity (see for instance \cite{BuBu}). 

Our point of view is to find results for rather general classes of functionals. In fact, the only conditions we require are the boundedness and the semicontinuity of the functionals. To obtain such a result and to successfully apply it to a large class of problems, we use the simplest possible topology on the selected class of domains, that is the so-called uniform convergence. This topology seems to be the ``most appropriate'' for the class of convex domains of $\R^N$ of prescribed volume. Moreover, it induces a rather strong convergence. By ``most appropriate'' we mean for instance that our abstract result can be adapted to a PDE of arbitrary order. 

This paper is inspired by the papers of G.~Buttazzo and B.~Kawohl and was realized with G.~Buttazzo after being awared of \cite{Ka}.

\section{Preliminary results}
In this section we first introduce some of the tools we shall need. In particular, we clarify or recall some properties of convex domains. 

We consider the following class of domains 
\[
C_m(D)=\{{\Om\subset D\ \mbox{open, convex and s.t. $|\Om|=m$}}\},
\]
where $D\subset \R^N$ is a compact set and  $m > 0$.
\begin{defi}\label{conv}
We say that a sequence $(\Om_n)_{n\in \N} \subset C_m(D)$ uniformly converges to 
$\Om \in C_m(D)$ if for every $\epsilon > 0$ one has
\begin{itemize}
\item[(d1)]
$\Om_n \subset \Om + B(\epsilon)$,
\item[(d2)] 
$\Om \subset \Om_n+B(\epsilon)$,
\end{itemize}
for $n$ large enough and with $B(\epsilon)$ the open ball of radius $\epsilon$ in $\R^N$ centered at the origin.
\end{defi}
In the following, a functional $J: C_m(D) \to \overline \R$ is said to be continuous if it is continuous according to this definition of convergence. 
\begin{rem}
Notice that the definition of the uniform convergence is equivalent to the convergence induced by the Hausdorff distance, 
\[
 d(\Om_1, \Om_2) = \max\{\sup_{x \in \Om_1}d(x,\Om_2),\sup_{x \in \Om_2}d(x,\Om_1)\}
\]
 with $\displaystyle d(x,A)=\inf_{y \in A}\{|x-y|\}$. It is also useful to define the complementary Hausdorff distance (see \cite{BuBu}), 
\[
 d_{H^c}(\Om_1, \Om_2) = \sup_{x \in \R^N}\ \bigl |d(x,\Om^c_1)-d(x,\Om^c_2) \bigr|.
\]
Similarly, for $(\Om_n)_{n \in \N} \subset C_m(D)$ and $\Om \subset C_m(D)$, we have, as $n \to \infty$, 
\[
\Om_n \to \Om \quad \mbox{uniformly}\quad \iff \quad d_{H^c}(\Om_n, \Om)\to 0.
\]
Moreover, it is well-known (see for instance \cite{Ku} or \cite{piron}) that the class of closed subsets of $D$ is compact for the Hausdorff convergence.
\end{rem}
We now give some properties of convex domains. First of all, convex domains have some kind of regular boundary. More precisely, according to \cite{Gri} (section 1.2) we have the following lemma.
\begin{lem}\label{pwsmooth} 
If $\Om$ is a convex domain $\subset \R^N$ then $\partial \Omega$ is locally
Lipschitz continuous. In particular, we have following properties of the boundary :
\begin{itemize}
\item[(i)] there exists a tangent plane ${\cal H}^{N-1}$-almost everywhere;
\item[(ii)] there are no zero interior angles.
\end{itemize}
\end{lem}
Note that points $(i)$ and $(ii)$ follow respectively from Rademacher's Theorem (see \cite{Ev}) and from the definition of convexity. According to \cite{Ev}, the boundary $\partial \Om$ can by described by a $W^{1,\infty}$ function of some parameters.

In addition to this property, any convex domain can be approximated by a sequence of  convex domains with smooth boundary.
\begin{lem}\label{smooth}
Let $\Om$ be a convex, bounded and open set of $\R^N$. Then for every $\epsilon >0$, there exist two convex open subset $\Om_1$ and $\Om_2$ in $R^N$ such that
\begin{itemize}
\item[(a)] $\Om_1 \subset \Om \subset \Om_2$,
\item[(b)] $\Om_1$ and $\Om_2$ have $C^2$ boundaries $\partial \Om_1$ and $\partial \Om_2$,
\item[(c)] $d(\partial \Om_1,\partial \Om_2) < \epsilon$, 
\end{itemize}
where $d(\partial \Om_1,\partial \Om_2)$ denotes the Hausdorff distance between
$\partial \Om_1$ and  $\partial \Om_2$. 
\end{lem}
{\bf{Proof.}} We refer to \cite{Gri}, lemma 3.2.1.1, p.147 for a proof.

\hfill {$\blacksquare$}

In the following lemmas, $|\Om|={\cal H}^N(\Om)$ denotes the ${\cal H}^N$-measure of $\Om$ and $|\partial \Om|={\cal H}^{N-1}(\partial \Om)$ denotes the ${\cal H}^{N-1}$-measure of its boundary (see \cite{Ambr} for an introduction to the Hausdorff measure).
\begin{lem}\label{convex}
The following properties hold for convex domains:
\begin{itemize}
\item[(i)] If $\Om_1 \subset \Om_2$ are two convex bodies then we also have $|\partial \Om_1 |\leq |\partial\Om_2|$.
\item[(ii)] If $(\Om_n)_n$ converges uniformly to $\Om$, then $|\Om_n| \to |\Om|$ and $|\partial \Om_n| \to |\partial \Om|$.
\end{itemize}
\end{lem}
{\bf{Proof.}} We refer to \cite{BuGu} for a proof.

\hfill {$\blacksquare$}

In the second statement, the convergence of the measures of the domains and of their boudaries is related to the convexity and to the Hausdorff convergence of the domains. However, the convexity is not a necessary condition for this property to hold. In fact, a more general class for this property would be the class of starshaped domains w.r.t a ball. Let us mention that in general, for open domains, if $d_{H^c}(\Om_n, \Om)\to 0$ as $n \to \infty$ then $\displaystyle \liminf_{n \to \infty}|\Om_n|\geq|\Om|$. Moreover, concerning the second part of statement $(ii)$, let us remark that for general open subsets of $\R^N$ even a relationship such as $\displaystyle \liminf_{n \to \infty}|\partial \Om_n|\geq|\partial \Om|$ is false in general. In fact, we only have the lower semicontinuity of the one-dimensional measure of the boundary w.r.t. the Hausdorff convergence of one-dimensional connected sets (see \cite{Ambr} for more results of this kind). 

Let us continue these preliminary results by a fundamental isoperimetric inequality, i.e. a relationship between volume and area of subsets of $R^N$. Such inequalities, called ``Bonnesen's Inequalities'' are summarized and proved in \cite {Oss}.
\begin{lem}\label{iso}
For convex bodies  $\Om \subset R^N, N \geq 2$, one has
\[
{\cal H}^N(\Om) < \rho \ {\cal H}^{N-1}(\partial \Om),
\]
where $\rho$ is the radius of the largest ball included in $\Om$. 
\end{lem}
{\bf{Proof.}} See \cite{Oss} for a proof.

\hfill {$\blacksquare$}

Let us compare this Bonnesen's Inequality to the following isoperimetric inequality (see \cite{EvGa})
\[
{\cal H}^N(\Om) < C \ \bigl |\bigl |\partial \Om\bigr |\bigr |^{\frac{N}{N-1}},
\]
valid in the more general context of sets of finite perimeter. Instead of the Haudorff $(N-1)$-dimensional measure on the right hand-side, one has the 'measure' perimeter. For convex domains we remark that the relationship ${\cal H}^{N-1}(\partial \Om)=\bigl |\bigl |\partial \Om\bigr |\bigr |$ holds (see \cite{BuGu}).                                                                                              
 
As announced, we seek to minimize the eigenvalues of some
elliptic operator $L$ of order $2p, p\geq 1$ over $C_m(D)$. For this reason we first need
to prove that every eigenvalue is a continuous functional on $C_m(D)$. To achieve
this aim, an important result is the well-known min-max
characterization of the $k$-th eigenvalue of $L$. In this paper, we mean by ellipticity of an operator that the associated weak form is coercive (sometimes this is called ``strong ellipticity'').
Let $\Om \in C_m(D)$ and let $H(\Om) \subset L^2(\Om)$ be a Hilbert space. If the
operator $L$ induces a symmetric and coercive weak form $(Lu,u)$ in $H(\Om)$, then the eigenvalues $\lambda_k \in \R$ and the eigenfunctions $u_k \in H(\Om)$ satisfy $Lu_k=\lambda_k u_k$ ($k \in \N$) with
\[
\lambda_k=\min_{H_k}\max_{\stackrel{u\in H_k,}{(u,u)=1}}(Lu,u),
\]
where  $H_k$ is any $k$-dimensional subspace of $H(\Om)$ and $(\cdot,\cdot)$ is the scalar product in $L^2(\Om)$.
Similarly, if the eigenvalues $\lambda_k$ and the eigenfunctions $u_k$
are such that $Lu_k=\lambda_k Bu_k$, for any elliptic operator $B$ of
order $2q$ with $q<p$ ($B$ is related to a symmetric and coercive weak form $(Bu,u)$),  we would have 
\[
\lambda_k=\min_{H_k}\max_{\stackrel{u\in H_k,}{(Bu,u)=1}}(Lu,u).
\]
Another important characterization reads as 
\begin{equation}\label{min}
\lambda_k=\min_{\stackrel{u\in H_{k-1},}{(Bu,u)=1}}(Lu,u)
\end{equation}
where $H_{k-1}$ is the subspace of $H$ defined with the ($k-1$) first eigenvectors $e_i$ as follows 
\[
H_{k-1}=\{{u \in H, (Bu,e_i)=0, 1\leq i \leq k-1}\}. 
\]
See \cite{RaTh},\cite{Wi} for a survey on this topic. This characterization will be very useful in proving the following monotonicity result. Let us define the space ${\cal A}(D)=\{{\Om\subset D,\ \Om \ \mbox{is open}}\}.$
\begin{lem}\label{monot}
If $(\lambda_k)_k$ is the sequence of all eigenvalues of the problem 
$Lu =\lambda Bu$, $u \in H^p_0(\Omega)$, for elliptic symmetric operators $L$ (of order $2p, p \geq 1$) and $B$ (of order $2q, q < p$), then the functional $\lambda_k :{\cal A}(D)\rightarrow \overline \R$ is monotone decreasing with respect to the set inclusion, i.e.
if $\Omega_1\subset \Omega_2$, then $\lambda_k(\Omega_1)\geq
\lambda_k(\Omega_2)$.
\end{lem}
{\bf Proof.} Since $H^p_0(\Omega)$ is the usual Sobolev space with all derivatives of order $\leq p-1$ being $0$ on the boundary $\partial \Om$, we remark that any competing function for $\lambda_k(\Omega_1)$
extended by $0$ in $\Omega_2\setminus \Omega_1$, is also a competing function
for $\lambda_k(\Omega_2)$. Then the proof follows from the characterization 
(\ref{min}) with $H=H^p_0$, since we have 
\[
\lambda_k(\Om_2)=\min_{\stackrel{v\in H_{k-1}(\Omega_2),}{(Bv,v)=1}}(Lv,v) \leq (L\tilde{u},\tilde{u}) = (Lu,u),
\]
for any $u\in H^p_0(\Omega_1)$ and its extension $\tilde{u}\in H^p_0(\Omega_2)$. When taking in the right hand-side the minimum over all normalized functions $u \in H_{k-1}(\Om_1)$, we obtain
\[
\lambda_k(\Om_2)=\min_{\stackrel{v\in H_{k-1}(\Om_2),}{(Bv,v)=1}}(Lv,v) \leq \lambda_k(\Om_1)=\min_{\stackrel{u\in H_{k-1}(\Om_1),}{(Bu,u)=1}}(Lu,u).
\]

{\hfill {$\blacksquare$}

For the existence of eigenvalues and eigenvectors to the spectral problem
$Lu=\lambda Bu$ we refer to spectral theorems, see for instance \cite{RaTh} or \cite{Wi} for a review.

This monotonicity result applies for instance to second order operators of the type  
\[
Lu=\frac{\partial}{\partial x_i}\left(a_{ij}(x)\frac{\partial u}{\partial x_j}\right)+ c_0(x)u,
\]
with symmetric and positive defined $a_{ij}(\cdot) \in {\cal C}(D)$ and nonnegative
$c_0 \in L^{\infty}(D)$.

\section{The Existence Result}\label{abstract}
The abstract existence result we shall use later is the following.
\begin{teo} \label{existence}
Let $D \subset \R^N$ be bounded and consider a functional
$J:C_m(D)\rightarrow \overline \R$ such that
\begin{itemize}
\item $J$ is bounded from below 
\item $J$ is lower semicontinuous , i.e. for each sequence $(\Om_n)_n$ such that $\Om_n \to \Om$ as $n \to \infty$,  we have $\displaystyle \liminf_{n\to\infty}J(\Om_n) \geq J(\Om)$. 
\end{itemize}
Then the infimum of  $J$  over  $C_m(D)$  is achieved, i.e. there exist  $\Om^{\star} \in C_m(D)$  such that 
\[
\inf_{C_m(D)}J(\Om) = \min_{C_m(D)}J(\Om)=J(\Om^{\star}).
\]
\end{teo}
{\bf Proof.} Since the functional is bounded from below, there exist a finite infimum in $C_m(D)$ and the result follows from the direct method of calculus of variations, as soon as we prove the set $C_m(D)$ is (sequentially) compact with respect to the uniform convergence. We then consider any minimizing sequence, and extract a subsequence which converges to the minimizer $\Om^{\star}$, the functional being lower semicontinuous.

Let us show the compactness of $C_m(D)$. If we consider any sequence $(\Om_n)_{n \in \N}$ in $C_m(D)$, then a ball of fixed radius $\rho >0$ fits every $\Om_n$, since, according to Lemma \ref{convex} and Lemma \ref{iso}, one has $\rho_n > \frac{|\Om_n|}{|\partial \Om_n|} > \frac{m}{S} > 0$, with $S$ the ${\cal H}^{N-1}$-measure of the boundary of $D$, and $\rho_n$ the largest inner radius associated to $\Om_n$. We define $\rho := \frac{m}{S}$ and consider the sequence of inner balls $(B_n)_{n \in \N}\subset D, B_n = B(x_n, \rho)$, with $(x_n)_{n \in \N}$ being the sequence of all the corresponding centers.
Since the sequence $(x_n)_{n \in \N}$ is contained in the compact set $D$, we extract from $(x_n)_n$ a subsequence converging to a point $x^{\star}$ of $D$.

Following an idea found in \cite{Cox}, we shall associate the above  sequence $(x_n)_{n \in \N}$ to the sequence of functions $(r_n)_n$, each radius $r_n$ describing parametrically the boundary $\partial \Om_n$ (this kind of parametric description make sense because the domains are convex).
Since $B(x_n,\rho) \subset  \Om_n \subset B[0,R_m]$ and $R_m$ is large enough such that $D \subset B[0,R_m]$, then by convexity of the $\Om_n$, the sequence $(\partial \Om_n)_{n \in \N}$ is equi-Lipschitz continuous, that is the Lipschitz constants of all $r_n$ are uniformly bounded, i.e. $L_n(x) < L$ for $x \in \partial \Om_n$. Denoting $S^{N-1}$ the unit ball of $\R^{N-1}$, the functions $r_n: S^{N-1}\rightarrow [\rho,2R_m]$ are such that $(r_n)_n$ is bounded in $W^{1,\infty}(S^{N-1})$. Hence we have by Rellich-Kondrachov's compact embedding Theorem that, up to a subsequence, $(r_n(\cdot))_n$  converges to a continuous function $r^{\star}(\cdot)$ in ${\cal C}(S^{N-1})$.
Since there is a sequence $((x_n,r_n(\cdot)))_n$ converging to $(x^{\star},r^{\star}(\cdot))$, then for given $\epsilon > 0$ we have, $n$ being large enough, $|x^{\star}-x_n|<\epsilon$ and $|r^{\star}-r_n|_{\infty}<\epsilon$. This is exactly, by definition \ref{conv} the uniform convergence. Thus ${\Om_n \to \Om^{\star}}$ uniformly, with $(\Om_n)_{n \in \N}$ a subsequence of the initial minimizing sequence and $B(x^{\star},\rho) \subset \Om^{\star} \subset D$. 
This convergence is strong enough to force $\Om^{\star}$ to be convex and according to Lemma \ref{convex} the volume is preserved, i.e. $|\Om^{\star}|=m$. Moreover, since $r^{\star}$ is continuous, its epigraph is closed and thus $\Om^{\star}$ is open. Since $\Om^{\star}$ is an admissible domain, this ends the proof.

\hfill {$\blacksquare$}

\begin{rem}
Clearly, the supremum of $J$ is achieved as soon as $J$ is bounded from above and upper semicontinuous, i.e. for each sequence $(\Om_n)_n$ such that $\Om_n \to \Om$ as $n \to \infty$,  we have $\displaystyle \limsup_{n\to\infty}J(\Om_n) \leq J(\Om)$. Moreover if the functional is bounded and continuous, then both the infimum and the maximum over $C_m(D)$ are achieved.
\end{rem}
\begin{rem}\label{volume}
Instead of the constraint $|\Om|=m$ we could prove the same existence result with the weaker condition that $c\leq |\Om|\leq m$, for any $c > 0$ and $m \leq |D|$.
\end{rem}
\begin{rem}\label{relax}
We could prove this result with other assumptions. According to the work of J.~Cox and M.~Ross (see \cite{Cox}), another suitable class of domains in $\R^2$ is that of starshaped domains which contain a disk, occupy a given area and do not exceed a prescribed parameter. There exists counter examples showing that these results for planar domains may not be extended to any space dimension $N >2$.
\end{rem}

\section{Minimum Problems for Eigenvalues.}

The monotonicity of the  eigenvalues (see Lemma \ref{monot}) is fundamental because it allows us to prove the continuity of the spectrum for a certain class of problems. For example, for an operator $L$ of order $2p$ with no lower order terms, the continuity follows from the monotonicity property and from the homogeneity of the eigenvalues, i.e. $\lambda_k(\alpha \Om)=\alpha^{-2p}\lambda_k(\Om)$ ($\alpha\in \R$). The homogeneity  holds if the $a_{ij}$ are constants. Unfortunately, the homogeneity fails at soon as the $a_{ij}$ are not constant functions but we will show in the following that the continuity of the eigenvalues still holds if the $a_{ij}(\cdot)$ are continuous functions on $D$. Moreover, we will mostly consider the case where the operator $B$ is the identity. 

Let us consider as first step the second order elliptic operator
of following type
\begin{equation}\label{deg2}
Lu=\frac{\partial}{\partial x_i}\left(a_{ij}\frac{\partial u}{\partial x_j}\right),
\end{equation}
where $(a_{ij})_{i,j}$ is a constant positive defined symmetric matrix of order $N \times N$. Since $\Om$ is bounded, then by Poincar\'e's Inequality the Sobolev space $H^1_0(\Om)$ is endowed with the norm $||u||_{H^1_0(\Omega)}^2\int_{\Om} u_{x_i}^2dx$ ($=\int_{\Om} (Du)^2dx$) and thus the operator $L$ admits a spectral decomposition. If $\Lambda(\Om)$ denotes the whole spectrum of $L$ over $H^1_0(\Om)$, we wish to prove that the problem
\begin{equation}\label{opt}
\min\{{\Phi(\Lambda(\Om)), \Om \in C_m(D)}\}
\end{equation}
admits at least a solution provided the functional $\Phi: \R^{\N}\to \overline \R$ is bounded from below and is lower semicontinuous in the sense that 
\[
\mbox{if} \quad \lambda_k(\Om_n)\to\lambda_k(\Om) \quad \mbox{for every $k$, as } n\to\infty \quad \mbox{then}  \quad \liminf_{n\to\infty}\Phi(\Lambda(\Om_n)) \geq \Phi(\Lambda(\Om)),
\]
where $\Lambda(\Om)=(\lambda_k(\Om))_{k \in \N}$ and $\Lambda(\Om_n)=(\lambda_k(\Om_n))_{k \in \N}$.
In order to prove this result, it is sufficient to prove the following lemma.
\begin{lem}\label{cont}
For an operator of type (\ref{deg2}) and for every $k \in \N$, the functional $\lambda_k: C_m(D)\to \overline \R$ is bounded and continuous.
\end{lem}
{\bf Proof.} The boundedness of $\lambda_k$ simply follows from the monotonicity of the eigenvalues (see Lemma \ref{monot}) because 
\[
\lambda_k(D)\leq \lambda_k(\Om) \leq \lambda_k(B_{\rho}) \qquad \forall \Om \in C_m(D),
\] 
where $B_{\rho}$ is an inner ball of radius $\rho$ (see proof of the Theorem \ref{existence}) included in $\Om$. We can consider, in this case, that $B_{\rho}$ is centered at the origin. To show the continuity we shall use the monotonicity and the homogeneity of the eigenvalues (since the $a_{ij}$ are constants, then by simply changing variable it is easy to show the homogeneity of degree $-2$). Since the domains are convex we have the existence of $t\geq 1$ such that for all $\epsilon >0$
\begin{equation}\label{inclusion}
\Om +B(0,\epsilon)\subset (1+t\epsilon)\Om,
\end{equation}
the expansion being performed with respect to $x_{\Om}\in \Om$ (for constant $a_{ij}$ we can consider $x_{\Om}$ to be the origin). Let us denote $\Om + B(0,\epsilon)$ by $\Om_{\epsilon}$.
By the relationship (\ref{inclusion}), by the homogeneity  and by Lemma \ref{monot} we have $(1+t\epsilon)^{-2}\lambda_k(\Om)=\lambda_k\bigl((1+t\epsilon)\Om\bigr) \leq \lambda_k(\Om_{\epsilon})$. By definition of the uniform convergence of the domains, for $n$ large enough $(d1)$ leads to $(1+t\epsilon)^{-2}\lambda_k(\Om) \leq \lambda_k(\Om_n)$ and $(d2)$ to $ (1+t\epsilon)^2\lambda_k(\Om) \geq \lambda_k(\Om_n)$. Hence, as $\epsilon \to 0$ we have proved the continuity of the $k$-th eigenvalue for every $k \in \N$.

\hfill {$\blacksquare$}

\begin{cor}\label{corollaire}
Problem (\ref{opt}) admits at least a solution.
\end{cor}
{\bf Proof.} The existence result follows now easily from Theorem \ref{existence}, since Lemma \ref{cont} allows us to show that $\Phi(\Lambda(\Om))$ is lower semicontinuous (l.s.c.) with respect to the uniform convergence as soon as $\Phi$ is  l.s.c. on $\R^{\N}$. 

\hfill {$\blacksquare$}

We consider now the Problem (\ref{opt}) associated with the spectral Problem (\ref{deg2}) where the $a_{ij}(\cdot)$ are not constants any more but are still continuous functions on $D$. In this case, we deal with second order operators of type
\begin{equation}\label{deg2nc}
Lu=\frac{\partial}{\partial x_i}\left(a_{ij}(x)\frac{\partial u}{\partial x_j}\right).
\end{equation}
We already pointed out that we do not have the homogeneity of the eigenvalues
\begin{equation}\label{exprmin}
\lambda_k(\Om)=\min_{\stackrel{u \in H_{k-1},}{\int_{\Om} u^2(x)dx=1}}\int_{\Om} a_{ij}(x)\frac{\partial u(x)}{\partial x_i}\frac{\partial u(x)}{\partial x_j}dx,
\end{equation}
with $H_{k-1}$ as defined in (\ref{min}) and the sum is taken over repeated indices. 

Our goal, is to prove an existence result to Problem (\ref{opt}) for an operator $L$ as (\ref{deg2nc}). The existence will easily follow from the following lemma.
\begin{lem}\label{contnc}
For every $k \in \N$, the functional $\lambda_k : C_m(D)\to\ {[-M,+M]}$ associated to the operator $L$ defined in (\ref{deg2nc}) is continuous, provided the $a_{ij}$ are continuous functions on D.
\end{lem}
{\bf Proof.} Since we have the characterization (\ref{exprmin}), $\lambda_n(\cdot)$ is monotone decreasing (see Lemma \ref{monot}) and the boundedness of $\lambda_k(\cdot)$ follows by the same arguments as in Lemma \ref{cont}. In fact, we have 
\[
\lambda_k(\Om)\leq\lambda_k\bigl(B(x_{\Om},\rho) \bigr)\leq \max_{i,j}\bigl | \bigl |a_{ij}\bigr | \bigr |_{L^{\infty}(D)}\hat \lambda_k \bigl(B(0,\rho)\bigr),
\]
with $\hat \lambda_k$ standing for the $k$-th eigenvalue of the problem with every  $a_{ij}(\cdot)$ set to 1.
To prove the continuity, let us fix $\epsilon >0$ and consider a sequence $(\Om_n)_{n\in \N}$ converging uniformly to $\Om$, then we have, by definition $(d2)$, by convexity of the domains and by equation (\ref{inclusion}) that for $n$ large enough, $(1+t\epsilon)^{-1}\Om \subset \Om_n$ with $t \geq 1$. Then
\[
\lambda_k(\Om_n) \leq \lambda_k\bigl((1+t\epsilon)^{-1}\Om\bigr)\qquad \mbox{(A)}.
\]
Since for $n$ large enough, $(1+t\epsilon)^{-1}\Om_n \subset \Om$ by definition $(d1)$ and relationship (\ref{inclusion}), it follows that
\[
\lambda_k(\Om)\leq \lambda_k\bigl((1+t\epsilon)^{-1}\Om_n\bigr) \qquad \mbox{(B)}.
\]
Moreover by (\ref{exprmin})
\[\lambda_k(G) \min_{\stackrel{u \in H_{k-1},}{\int_{G} u^2(x)dx=1}}\int_{G}
a_{ij}(x)\frac{\partial u(x)}{\partial x_i}\frac{\partial
u(x)}{\partial x_j}dx=\int_{G}
a_{ij}(x)\frac{\partial u_G(x)}{\partial x_i}\frac{\partial
u_G(x)}{\partial x_j}dx,
\]
with the generic symbol $G$ used for both
$\Om_n$ and $\Om$ and $u_G\in H_{k-1}$, a solution of $Lu(x)=\lambda_k(G)
u(x)$ for $x \in G$. Hence, dealing first with $\Om$ we have 
\[
\lambda_k\bigl((1+t\epsilon)^{-1}\Om\bigr) \leq \frac{\displaystyle \int_{(1+t\epsilon)^{-1}\Om} a_{ij}(x)\frac{\partial
u^{\epsilon}_{\Om}(x)}{\partial x_i}\frac{\partial
u^{\epsilon}_{\Om}(x)}{\partial x_j}dx}{\displaystyle \int_{(1+t\epsilon)^{-1}\Om}
|u^{\epsilon}_{\Om}(x)|^2dx} ,
\]
with the function $u^{\epsilon}_{\Om}(\cdot) = u_{\Om}\bigl((1+t\epsilon)\cdot\bigr)$ being a competing function for the infimum on $(1+t\epsilon)^{-1}\Om$. Then, by a mere change of  variables in the right hand-side, it follows that 
\begin{equation}\label{ineg}
\lambda_k\bigl((1+t\epsilon)^{-1}\Om\bigr) \leq (1+t\epsilon)^{2}\frac{\displaystyle \int_{\Om} a_{ij}(\frac{x}{1+t\epsilon})\frac{\partial
u_{\Om}(x)}{\partial x_i}\frac{\partial
u_{\Om}(x)}{\partial x_j}dx}{\displaystyle \int_{\Om}
|u_{\Om}(x)|^2dx}
\end{equation}
By continuity of the $a_{ij}(\cdot)$, it then follows that
\[
\mbox{(C)} \qquad \liminf_{\epsilon\to 0}\lambda_k\bigl((1+t\epsilon)^{-1}\Om\bigr) \leq \lambda_k(\Om) \qquad \mbox{and by (A)}\qquad \liminf_{n\to \infty}\lambda_k(\Om_n) \leq \lambda_k(\Om)\qquad \mbox{(D)}.
\] 
Inequality (C) holds with $\Om_n$ instead of $\Om$, so from(B) we have
\[
\lambda_k(\Om) \leq \liminf_{n\to\infty} \liminf_{\epsilon \to 0} \lambda_k((1+t\epsilon)^{-1}\Om_n)\leq \liminf_{n\to\infty}\lambda_k(\Om_n). 
\]
It follows then from this last inequality and from (D) that 
\[
\liminf_{n\to \infty}\lambda_k(\Om_n) =\lambda_k(\Om)= \limsup_{n\to \infty}\lambda_k(\Om_n)= \lim_{n\to \infty}\lambda_k(\Om_n),
\]
since the $\displaystyle \limsup$ could be used instead of the $\displaystyle \liminf$ for the above limit process. The continuity is proved.

\hfill {$\blacksquare$}

\begin{cor}\label{cor4.4}
Problem (\ref{opt}) associated to the spectral decomposition of (\ref{deg2nc}) admits at least a solution.
\end{cor}
{\bf Proof.} See the proof of Corollary \ref{corollaire}.

\hfill {$\blacksquare$}

Our existence result allows us to consider any symmetric (strongly) elliptic operator of order $2p, p\geq 1$, since the uniform convergence is a priori adapted to any Sobolev space $H_0^p(\Om), p \geq 1$.
For example, let us consider a 4-th order operator $L$ of the form 
\begin{equation}\label{deg4}
Lu=\frac{\partial^2}{\partial x_i\partial
  x_j}\left(a_{ijkl}(x)\frac{\partial^2u}{\partial x_k\partial x_l}\right).
\end{equation} 
with $a_{ijkl}(\cdot)$ a continuous tensor on $D$ verifying the symmetry assumption $a_{ijkl}=a_{ijlk}=a_{jikl}=a_{klij}$ and such that $L$ is (strongly) elliptic, i.e. there exists $\nu>0$
such that for every $\xi \in \R^{N \times N}$
\[
\nu \xi_{ij}^2\leq a_{ijkl}(x)\xi_{ij}\xi_{kl} \qquad \mbox{a.e. in $D$}.
\]
Since $\Om$ is bounded, by Poincar\'e's Inequality the Sobolev space $H^2_0(\Om)$ is endowed with the norm $||u||_{H^2_0(\Omega)}^2\int_{\Om}u_{x_ix_j}^2dx$ ($=\int_{\Om}|D^2u|^2dx$) whence the operator $L$ admits a spectral decomposition in $H^2_0(\Omega)$.

The proof of a solution for Problem (\ref{opt}) relies again on the continuity of every eigenvalue
\begin{equation}
\gamma_k(\Om)=\min_{\stackrel{u \in H_{k-1},}{\int_{\Om} u^2(x)dx=1}}\int_{\Om}a_{ijkl}(x)\frac{\partial^2u}{\partial x_i\partial x_j}\frac{\partial^2u}{\partial x_k\partial x_l}dx.
\end{equation}
This can be easily proved as in Lemma \ref{cont} or Lemma \ref{contnc}. 

We point out that the method used for proving Corollary \ref{corollaire} and Corollary \ref{cor4.4} can be used in a similar way to easily prove such existence results for any linear, symmetric and strongly elliptic operators $L$ of order $2p$, provided we have the continuity of the coefficients of $L$ on $D$ and seek an optimal domain over the class of all bounded convex domains with prescribed measure. 
\begin{rem}
The continuity property also holds for a class of operators such as
(\ref{deg2}) or (\ref{deg4}) with an additional  zero order coefficient $c_0(x)u$, provided $c_0$ is nonnegative and is a continuous function on $D$. In fact, for a second order operator,  relation (\ref{ineg}) becomes 
\[
\lambda_k((1+t\epsilon)^{-1}\Om) \leq \qquad \mbox{...} \qquad
+ \frac{\displaystyle \int_{\Om} c_0(x)(u_{\Om}(x))^2dx}{\int_{\Om}
|u_{\Om}(x)|^2dx}
\] 
and the result follows easily.

Finally, the same kind of arguments can also be repeated for a larger class of operators, having the general form
\begin{equation}\label{compl}
\frac{\partial^2}{\partial x_i\partial
  x_j}\left(a_{ijkl}(x)\frac{\partial^2u}{\partial x_k\partial x_l}\right)-
\frac{\partial}{\partial x_i}\left(b_{ij}(x)\frac{\partial u}{\partial x_j}
\right)+c_0(x)u,
\end{equation}
with nonnegative, continuous and symmetric $b_{ij}$ and nonnegative, 
continuous $c_0$. If for every $(i,j)$, $b_{ij}$ and
$c_0$ are not nonnegative on $D$, then these low order terms may cause problems concerning the coerciveness of the weak form, unless
$|b_{ij}|$, $|c_0|$ and the measure of $\Omega$ are ``small enough''
compared to  the ellipticity constant $\nu$ (see for instance \cite{La}).
\end{rem}
\begin{rem}
The following is known (see \cite{BuBu}) when the domains are not required to be convex but are all quasi-open sets included in $D$. For $L=-\Delta$, any lower
semicontinuous functional of the variable $\bigl( \lambda_1(\Om),\lambda_2(\Om) \bigr)$ achieve its infimum. Moreover, if the functional is also monotone decreasing with respect to the set inclusion, then this last result also holds for the whole spectrum of $-\Delta$.
\end{rem}
\begin{rem}\label{dimi}
We considered the first application for the general Sobolev spaces $H^p_0(\Om)$. For an operator $L$ such as \ref{deg4} in the space $H^1_0(\Om) \cap H^2(\Om)$, the previous results do not hold, since the monotonicity property of the eigenvalues fails.
\end{rem}

\section{Minimizing Functionals given by an Integral.}

In this second application let us consider functionals of following type
\[
J: C_m(D) \to \R \qquad \mbox{is such that} \qquad J(\Om) = \int_Dj\bigl(x,u_{\Om},Du_{\Om}(x),...,D^{p-1}u_{\Om}(x)\bigr)dx,
\] with $j: \Om \times \R^{p-1} \to \overline \R$ and where $u_{\Om}$ is the unique solution, extended by zero outside $\Om$, of the problem
\[
u \in H^p_0(\Om), \qquad Lu = f, \qquad f \in L^2(D),
\]
with the operator $L$ of order $2p, p \geq 1$. Let us denote the solution as $u_{\Om}=L^{-1}(\Om,f)$.
We will consider the problem 
\begin{equation}\label{optL}
\min_{\stackrel{\Om \in C_m(D)}{u_{\Om}=L^{-1}(\Om,f)}} \int_D j\bigl(x,u_{\Om}(x),Du_{\Om}(x),...,D^{p-1}u_{\Om}(x)\bigr)dx.
\end{equation}
We are interested in searching for a solution to the Problem (\ref{optL}) using the  abstract Theorem \ref{existence}. We remark that without the convexity constraint on the domains, problem (\ref{optL}) could have its infimum achieved outside the set we have chosen for the minimizers. In this case, one should relax the problem, i.e. minimize a generalized functional over more general sets, for instance a class of measures, and try to link this latter problem to the original one, as done in \cite{BuDa} for a class of  open but not necessarily convex domains of $\R^N$ and for $p=1$. We shall nevertheless prove thanks to convexity of admissible domains that Problem (\ref{optL}) admits at least a solution provided the function $j: \Om \times \R^{p-1} \to \R$ is bounded from below and lower semicontinuous in the $(p-1)$ last variables. We remark that no growth condition on $j$ is required in our assumptions. We shall consider second order operators of type (\ref{deg2nc}) and fourth order operators of type (\ref{deg4}) and finally apply the existence result to the biharmonic operator $\Delta^2$. Let us start with $L$ as in (\ref{deg4}). For more general operators $L$ of higher order, the proof is similar to the following one, but for convenience we restrict ourselves to the fourth-order situation.

\begin{cor}\label{born1}
The optimization Problem (\ref{optL}) with the 4-th order operator $L$ as in 
(\ref{deg4}) admits a solution.
\end{cor}

{\bf Proof.} The weak form of the equation $Lu=f$ on $\Om$ is, for $u \in H^2_0(\Om)$
\begin{equation}\label{weak}
(Lu,\phi)=\int_{\Om} a_{ijkl}(x)\frac{\partial^2u(x)}{\partial x_i\partial x_j}\frac{\partial^2\phi(x)}{\partial x_k\partial x_l}dx = 
\int_{\Om} f(x) \phi(x) dx \qquad \forall\phi \in C^{\infty}_0(\Om)
\end{equation}
The tensor $a_{ijkl}(\cdot)$ is continuous on the compact $D$ and thus it is also bounded. Since it is also positive definite on $D$, the bilinear and symmetric form $(Lu,v)$ defines a scalar product on $H^2_0(\Om)$, i.e. 
\begin{equation}\label{ps}
\nu\int_{\Om}u^2_{x_i x_j}(x)dx \leq (Lu,u) \leq \alpha \int_{\Om}u^2_{x_i x_j}(x)dx, 
\end{equation}
with $\nu,\alpha>0$. The existence of a unique weak solution $u_{\Om}$ in $H^2_0(\Om) $ for the equation $Lu=f$ on $\Om$ follows then by Lax-Milgram's Theorem. Problem (\ref{optL}) is well defined and we chose a minimizing sequence $(\Om_n)_{n \in \N}$.
Let us show now the solutions $u_{\Om_n}\to u_{\Om}$ in $H^1_0(D)$ provided $\Om_n$ uniformly converges to $\Om$. Since we have (\ref{weak}) and (\ref{ps}), we also have
\[
\nu ||u_{\Om_n}||^2_{H^2_0(D)}=\nu\int_D u^2_{x_i x_j}(x)dx \leq(Lu_{\Om_n},u_{\Om_n}) = (f, u_{\Om_n})
\]
\[
\leq ||f||_{L^2(\Om_n)}||u_{\Om_n}||_{L^2(\Om_n)}\leq ||f||_{L^2(D)}||u_{\Om_n}||_{H^2_0(D)}
\]
and finally $||u_{\Om_n}||_{H^2_0(D)}\leq \frac{||f||_{L^2(D)}}{\nu}$ with the solutions $u_{\Om_n}$ set equal to 0 outside $\Om_n$. Thus
$u_{\Om_n}\to u$ weakly in $H^2_0(D)$ and by Rellich's Theorem, $u_{\Om_n}\to u$ strongly in $H^1_0(D)$. It remains to show $u=u_{\Om}$. Since the domains are convex, the uniform convergence of the domains implies that any test function $\phi \in {\cal C}^{\infty}_0(\Om)$ is a test function in ${\cal C}^{\infty}_0(\Om_n)$, provided $n$ is large enough. Hence, for $\phi \in {\cal C}^{\infty}_0(\Om)$ and $n$ large enough, 
\[
(f,\phi)=(Lu_{\Om_n},\phi)=\int_D a_{ijkl}(x)\frac{\partial^2u_{\Om_n}(x)}{\partial x_i\partial x_j}\frac{\partial^2\phi(x)}{\partial x_k\partial x_l}dx, 
\]
by the above property of the test function. This integral converges to
\[
\int_D a_{ijkl}(x)\frac{\partial^2u(x)}{\partial x_i\partial x_j}\frac{\partial^2\phi(x)}{\partial x_k\partial x_l}dx=(Lu,\phi)
\]
as $n \to \infty$, by the weak convergence of the $u_{\Om_n}$ in $H^2_0(D)$ . Hence
\[
(Lu,\phi)=(f,\phi)\qquad\forall\phi\in C^{\infty}_0(\Om).
\]
Moreover, we have the ${\cal H}^N$-almost everywhere pointwise convergence of $u_{\Om_n}(x)$ to $u(x)$ and of $(Du)_{\Om_n}(x)$ to $(Du)(x)$ as $n \to 0$, that is, $u(x)=(Du)(x)=0$ ${\cal H}^N$-a.e. on $\partial \Om$ by the classical trace theorems. Since the weak solution of $Lu=f$ is unique in $H^2_0(\Om)$, it follows that $u=u_{\Om}$.

Thus, given $f\in L^2(D)$, we proved that both $u_{\Om_n} \to u_{\Om}$ and $(Du)_{\Om_n} \to (Du)_{\Om}$ in $L^2(D)$.  
It follows that for almost every $x \in D$ there are two subsequences $u_{\Om_n}(x) \to u_{\Om}(x)$ and $(Du)_{\Om_n}(x) \to (Du)_{\Om}(x)$. By Fatou's Lemma and by the semicontinuity of $j$,
\[
\liminf_{n \to \infty}\int_{D}j\bigl(x,u_{\Om_n}(x),Du_{\Om_n}(x)\bigr)dx \geq \int_{D}\liminf_{n \to \infty}j\bigl(x,u_{\Om_n}(x),Du_{\Om_n}(x)\bigr)dx 
\]
\[
\geq \int_{D}j\bigl(x,u_{\Om}(x),Du_{\Om}(x)\bigr)dx,
\]
whence
\[
\liminf_{n \to \infty}J(\Om_n) \geq J(\Om). 
\]
The functional being bounded from below, the existence result follows from the lower semicontinuity of the functional and from the existence Theorem \ref{existence}.

\hfill {$\blacksquare$}

Again, we easily obtain a generalization of Corollary \ref{born1} in the case that the operator $L$ is of order $2p$, is linear, symmetric and strongly elliptic. In addition to these assumptions,  required conditions for the proof are the boundedness of the coefficient of $L$ on $D$ and the convexity of the domains. Note that the latter condition may be weakened (see for instance the definition of the ``compactivorous property'' for connected domains in the paper \cite{BuZo}). 

\begin{rem}
If $j(x,0,\cdot \cdot \cdot,0)=0$ a.e. in $D$, then the previous existence result holds for 
$\displaystyle J(\Om) = \int_{\Om}j\bigl(x,u_{\Om},Du_{\Om}(x),...,D^{p-1}_{\Om}(x)\bigr)dx$ since the solutions are as usual set equal to $0$ outside $\Om$.
\end{rem}

\begin{rem}
Corollary \ref{born1} still holds for a generalized second order operator $L$ with first order and zero order terms (such as (\ref{compl})). In fact, we can follow the proof of Corollary \ref{born1} provided the low order coefficients
and the size of $\Om$ are ``small enough'' compared to  the ellipticity
constant $\nu$. 
\end{rem}

Let us apply Corollary \ref{born1} to the biharmonic operator, i.e let us consider optimization problem (\ref{optL}) associated to the equation
\[
\left\{\begin{array}{ll}
\Delta^2u=f  & \mbox{in }\Omega,\\
u=\frac{\partial u}{\partial \nu}=0 & \mbox{on }\partial \Omega
\end{array}\right.
\]
with $\Om \in C_m(D)$ and $f \in L^2(D)$. 

It is enough to show that the operator $\Delta^2$ verifies (\ref{ps}) (i.e. that $\Delta^2$ is strongly elliptic) in Corollary \ref{born1} (here we have $a_{ijkl}=\delta_{ij}\delta_{kl}$). In fact, in \cite{La} the authors proved (easily) that $||u||_{H^2_0(\Om)} \leq \beta ||\Delta u||_{L^2(\Om)} \leq \beta ||u||_{H^2_0(\Om)}$ and Corollary \ref{born1} then follows for the biharmonic operator. 

\begin{rem}
We could wonder whether the second application holds for spectral problems, i.e. for 
\[
u \in H^p_0(\Om), \qquad Lu = \lambda u, \qquad \lambda \in \R,
\]
with the operator $L$ of order $2p, p \geq 1$. The question is once again the existence of a minimizer of the functional $J(\Om) = \int_Dj\bigl(x,u^k_{\Om},Du^k_{\Om}(x),...,D^{p-1}u^k_{\Om}(x)\bigr)dx$ for the $k$-th eigenfunction $u^k$.
Actually, we could follow the previous proof if we knew the $k$-th eigenvalue $\lambda_k(\Om)$ is single, or equivalently if the associated eigenfunction $u^k_{\Om}$ is unique modulo a constant (for a strongly elliptic and symmetric $L$, the $k$-th eigenfunction on $H^p_0(\Om_n)$ is such that $||u_k||^2_{H_0^p(\Om_n)}\leq\frac{\lambda_k(\Om_n)}{\nu}\leq \frac{\lambda_k\bigl(B(\rho)\bigr)}{\nu}$, for $n$ large enough). Thus the previous existence result holds for instance for $L=\Delta$ in $H^1_0(\Om)$ and $k=1$. For other problems, for instance for the Laplacian and $k \geq 2$, for the Neumann problem involving the Laplacian, for the bilaplacian in $H^2_0(\Om)$ or for the Lam\'e system of linear elasticity, the existence of a minimizer still holds for the following problem
\[
\min_{\stackrel{\Om \in C_m(D)}{u^k_{\Om} \in L^{-1}(\Om,k)}} \int_D j\bigl(x,u^k_{\Om}(x),Du^k_{\Om}(x),...,D^{p-1}u^k_{\Om}(x)\bigr)dx,
\]
with $L^{-1}(\Om,k)$ being a finite dimensional subset of $H^p_0(\Om)$.
\end{rem}

\section{Newton's Problem of Minimal Resistance.}

We are interested in proving the existence of a convex body $\Om \subset \R^3$ which minimizes the so-called total resistance (see \cite{BuKa} and \cite{Ka} and the references therein). According to Newton's model, the total resistance is given by
\[
R(G,u)=\int_G \frac{1}{1+|Du|^2}dx, 
\]
with $u: G \to \R$ and for any given domain $G \subset \R^{2}$.
For a convex $G$ and a concave $u$, a $3$-dimensional convex body $\Om$ can be built with $G$ and the graph of $u$. In \cite{BuKa} the authors proved that the functional $R$ achieves its infimum over the class $C_M=\{u$ concave on $G: 0 \leq u \leq M\}$. 

In \cite{BuGu} the authors introduced another formalism for writing the total resistance in $\R^N$. They showed that $R(G,u)$ could be re-written as
\begin{equation}\label{resist}
F(\Om)=\int_{\partial \Om}((\nu\cdot A)^{+})^3d{\cal H}^{N-1},
\end{equation}
with $(\nu\cdot A)^{+}$ denoting the positive part of the projection of the normal $\nu$ on the stream direction $A$. Notice that we consider only the positive part of the scalar product, i.e. $u^{+}:=\frac{1}{2}\bigl ( |u|+u\bigr )$ in order to eliminate the points that are not relevant for the total resistance. With this notations, the authors proved the existence of a minimizer of the functional $F$ on the class
\[
C'_m(D)=\{ \Om \subset D, \mbox{open and convex and s.t. } m \leq |\Om|\}.
\]
For proving this result, they made use of tools from geometric measure theory, starting with an abstract theorem of Reshetnyak. Our goal is to prove this result using the tools and the results developed above. Once again, the convexity will be the most relevant property of the problem. 
\begin{prop}\label{contnorm}
Let $\Om_n$ be a sequence of convex domains converging to $\Om$ uniformly. Then for all $x \in \partial \Om \setminus S$ and all $x_n \in \partial \Om_n \setminus S_n$ we have
\[
\mbox{if} \qquad x_n \to x  \qquad \mbox{then} \qquad \nu_n(x_n) \to \nu(x), 
\]
where $S_n$ and $S$ respectively denote the sets of points where $\partial \Om_n$ and $\partial \Om$ are not differentiable.
\end{prop} 

{\bf{Proof.}} Since $\Om_n \to \Om$ as $n \to \infty$, then for $n \geq n_0$ we have a ball $B(x,\rho)$ included in every $\Om_n, n \geq n_0$. Let us consider a point $x^{\star} \in \partial \Om$ where the boundary is differentiable, i.e. where the tangent plane $T(x^{\star})$  exists. By Lemma \ref{pwsmooth} there is a tangent plane almost everywhere, i.e. the sets $S_n$ and $S$ are ${\cal H}^{N-1}$ - negligible. Since all bodies are convex, the half line starting at $x$ and passing by $x^{\star}$ crosses each $\Om_n$ at only one point $x_n$. If we discard all $x_n$ where the boundary $\partial \Om_n$ has a corner, it remains a sequence $(\hat x_n)_{n \geq n_0}$ of points of differentiability converging to $x^{\star}$ as $n \to \infty$ by definition of the uniform convergence. If we can not extract from $(x_n)_{n \geq n_0}$  a subsequence of points of differentiability converging to the point $x^{\star}$, then we discard the point $x^{\star}$. The points we discarded remain anyway ${\cal H}^{N-1}$ - negligible. In the following, the sets $\hat \Om_n$ and $\partial \hat \Om_n$, etc, are related to the sequence $(\hat x_n)_{n \in \N}$.

Let us define the tangent plane $\hat T_n(\hat x_n)$ to $\partial \hat \Om_n$ at $\hat x_n$. We have defined a sequence of points of differentiability which all lie on a straight line, then by passing to a subsequence if necessary, we can consider either they all lie between $x$ and $x^{\star}$ or all above $x^{\star}$. We first consider the second  situation occurs. 
If $\hat r_n$ is the parametric description of the $\partial \hat \Om_n$ then, since $B(x,\rho) \subset \hat \Om_n \subset D$, we know the sequence $(\nabla \hat r_n(\cdot))_{n \geq n_0}$ is equi-bounded and belongs then to a bounded subspace of $L^p(S^{N-1})$, for some $1 \leq p < \infty$. Moreover, for $n \geq n_0$ and for given $\epsilon > 0$, then for every small enough $h \geq 0$ we have by Lebesgue's theorem that ${|\nabla \hat r_n(\cdot + h)-\nabla \hat r_n(\cdot)}|_{L^p} \leq \epsilon$. Then, by Riesz-Kolmogorov's compactness theorem  the existence of a converging subsequence in $L^p(S^{N-1})$ follows (see \cite{Bre}) easily. Let us write $\nabla \hat r_n(\cdot) \to \gamma(\cdot)$ in $L^p(S^{N-1})$. Then, up to a subsequence, $(\nabla \hat r_n(\cdot))_{n \geq n_0}$ converges almost everywhere on $S^{N-1}$, i.e. $\nabla \hat r_n(\hat x_n) \to \gamma(x^{\star})$ as $n \to \infty$ or, in terms of normal vectors, $\hat \nu_n(\hat x_n) \to \hat \nu(x^{\star})$ as $n \to \infty$. It remains to prove that $\hat \nu(x^{\star})=\nu(x^{\star})$, where $\nu(x^{\star})$ is the normal vector to $\Om$ at $x^{\star}$. 
Indeed, assume $\hat \nu(x^{\star})\not =\nu(x^{\star})$, then the normal plane to $\hat \nu(x^{\star})$ at $x^{\star}$ cuts the domain $\Om$  in two parts and lies below a portion of the boundary $\partial \Om$ of area $\delta >0$. This creates a contradiction. In fact, since the converging domains $\hat \Om_n$ are convex, the tangent plane at $\hat x_n$ (orthonormal to $\hat \nu_n$) lies above the boundary $\partial \hat \Om_n$. Then, by the uniform convergence of the domains, the boundaries $\partial \hat \Om_n$ converge to $\partial \Om$ uniformly, and at the same time the tangent planes $\hat T_n$ also converge uniformly to $T$, with $T$ the normal plane to $\nu$ at $x^{\star}$. This would mean that $\delta \to 0$, which is a contradiction.  The same kind of arguments holds if the sequence of points $(\hat x_n)_{n \geq n_0}$ lie above $x^{\star}$.
This ends the proof.

\hfill {$\blacksquare$}

With this result it is now easy to prove the main result of this section.
\begin{prop}\label{derniere}
Among  all convex bodies of $D \subset \R^N$ of prescribed volume, there exists at least one body for which 
\[
F(\Om)=\int_{\partial \Om}f(x,\nu(x))d{\cal H}^{N-1}
\]
is minimized, provided the function $f$ is measurable and lower semicontinuous in the second variable.
\end{prop}

{\bf{Proof.}} Since for every $n$, there exists a one-to-one application $\phi_n: \partial \Om^{\star} \to \partial \Om_n$ almost everywhere differentiable and s.t. $\phi_n^{-1}$ is a.e. differentiable, we have
\[
F(\Om_n)=\int_{\partial \Om_n}f(x_n,\nu_n(x_n))d{\cal H}^{N-1}=\int_{\partial \Om^{\star}}f(x_n,\nu_n(x_n))|J_{\phi_n}|d{\cal H}^{N-1}.
\]
Then, according to Proposition \ref{contnorm} and by Fatou's Lemma, the functional $F:C'_m(D) \to \overline \R$ is lower semicontinuous with respect to the uniform convergence. Thus, the proof follows from the existence result of Section \ref{abstract} and from Remark \ref{volume}.

\hfill {$\blacksquare$}

\begin{cor}
Among all convex bodies of $D \subset \R^N$ of prescribed volume, there exists at least one body for which the total resistance (\ref{resist}) is minimum.
\end{cor}
{\bf{Proof.}} We apply Proposition \ref{derniere} to $f(x,\nu(x))=\bigl((\nu(x)\cdot A(x))^{+}\bigr)^3$ with the unit vector field $A$ being the stream direction.

\hfill {$\blacksquare$}

\begin{ack}
The author is very grateful to G.~Buttazzo who proposed and supervised this work. He would like to express his appreciation of the hospitality of the Dipartimento di Matematica dell'  Universit\`a di Pisa. He was supported by a CGRI grant (Communaut\'e Fran\c caise de Belgique) during the period from October 1999 until August 2000.
\end{ack}
\small author's e-mail adress: $vangoeth@mema.ucl.ac.be$. Visit also $http://cvgmt.sns.it$.
\ifx\undefined\allcaps\def\allcaps#1{#1}\fi

\end{document}